\documentclass[12pt]{article}
\usepackage{amsmath, amsfonts}
\usepackage{amssymb}
\usepackage{mathrsfs}

\advance\textheight by 5\baselineskip

\date{}

\begin{document} 

\begin{center}
    {\large \bf Projections of orbital measures for classical Lie groups\footnote{This work is supported by the Russian Science Foundation under grant 14-50-00005.}}
    
    {Dmitry Zubov}
\end{center}

{\bf 1. Introduction.}
    Let $G_n$ be a compact Lie group from the list: $ \mathrm{SO}(2n+1)$, $\mathrm{Sp}(2n), $ $\mathrm{O}(2n)$. Consider the adjoint action of $G_n$ on its Lie algebra $\mathfrak{g}_n$. Since $G_n$ is a compact group, for each orbit of the action there exists a unique $G_n$-invariant probability measure supported on this orbit. We shall call this measure {\it the orbital measure}.  

    Each orbit can be parametrized by a $n$-tuple $X = (x_1 \le ... \le x_n)$, $x_1 \ge 0$ of weakly decreasing numbers, which corresponds to the canonical form of the matrix from $\mathfrak g_n$.
    Let us denote by $\mathcal{X}_n$ the set of such $n$-tuples.
    
    Now we consider a natural projection map $p_k^n: \mathfrak{g}_n \to \mathfrak{g}_k$. Let $\mu_X$ be a $G_n$-orbital measure; then the measure $p_k^n(\mu_X)$ is invariant under the action of the group $G_k$. Each invariant measure can be represented as a continuous combination of orbital measures: in fact, for each Borel subset $S \in \mathfrak{g}_k$ we have
    $$
	p_k^n(\mu_X) (S) = \int\limits_{Y \in \mathcal{X}_k} \mu^{(k)}_Y(S) \cdot \nu_{X,k}(dY),
$$ 
    where $\mu^{(k)}_Y,  Y \in \mathcal{X}_k$ are $G_k$-orbital measures, while $\nu_{X,k}$ is a certain probability measure on the set $\nu_{X,k}$ of $k$-tuples.
    
    The measure $\nu_{X,k}$ is called {\it the radial part} of the measure $p_k^n(\mu_X)$. In the case of the unitary group $\mathrm{U}(n)$ the measure $\nu_{X,k}$ was computed by Olshanski \cite{Olsh} and Faraut \cite{Faraut}. Using the method of Faraut, we compute these measures for the groups $\mathrm{SO}(2n+1)$, $\mathrm{Sp}(2n), $ $\mathrm{O}(2n)$, expressing the radial parts of the orbital measures in terms of determinants of B-splines with the knots that are symmetric with respect to 0.
    
    \medskip
    
{\bf 2. Main result}
Before we formulate the main result, let us give the necessary definitions.

According to Curry and Schoenberg, \cite{CSh}, the {\it B-spline} with $n$ knots $t_1 < ... < t_n$ is a  ${\mathit C}^{n-3}$-smooth function $M_n(t_1,...,t_n; t)$ on $\mathbb R$ with the following properties: \par
1) $\mathrm{supp}(M_n(t_1,...,t_n; t)=[t_1,t_n]$;\par
2) The function $M_n(t_1,...,t_n; t)$ is the polynomial in $t$ of degree $n-2$ on each subinterval $(t_i, t_{i+1})$;\par
3) $\int_{\mathbb{R}} M_n(t_1,...,t_n;t) dt = 1$.\par

Note that the conditions above define the B-spline $M_n(t_1,...,t_n; t)$ in a unique way.

Recall also that {\it the divided differences} of a function $f$ are defined by induction as follows:
$$
	f[t_1,t_2] = \frac{f(t_1)-f(t_2)}{t_1 - t_2};   \, ...; \, f[t_1,...,t_n] = \frac{f[t_1,t_2,...,t_{n-1}]-f[t_2,...,t_{n-1},t_n]}{t_1-t_n}.
	$$
    
    The Hermite-Genocchi formula connects B-splines with divided differences. This formula tells that for a function $f$ which has a piecewise continuous derivative of degree $n-1$ there is an equality
    $$
        f[t_1,...,t_n] = \frac{1}{(n-1)!}\int_\mathbb{R} f^{(n-1)}(t)M_n(t_1,...,t_n; t) dt.
    $$
    
    Using the method of Faraut \cite{Faraut},for the groups $\mathrm{SO}(2n+1), \mathrm{Sp}(2n), \mathrm{SO}(2n)$  we obtain the formulas for $\nu_{X,k}$ in terms of determinants of B-splines (Theorem 1).

    Let us denote by
    $$V_n(T) := V_n(t_1,...,t_n) = \prod\limits_{1\le i<j \le m} (t_i-t_j)$$
    the Vandermonde determinant in the variables $T=(t_1,...,t_n)$, and let us also set
    $$
	v_n(dT) = V_n(T^2) \cdot \prod_{i=1}^n dt_i = V_n(t_1^2,...,t_n^2) \cdot \prod_{i=1}^n dt_i.
    $$
    Finally, for arbitrary $j,m \in \mathbb N$ we shall use a shortened notation
    $$
	M_{2m+2}(\pm t|_j^{m+j};t) := M_{2m+2}(-t_{m+j},...,-t_j,t_j,...,t_{m+j}; t).
$$

{\bf Theorem 1.}
    {\it The radial part $\nu_{X,k}$ of projection of the $G_n$-orbital measure $\mu_X$ is given by the following formula:
    $$
	\nu_{X,k}(dY) = \frac{c(n,k)}{\prod\limits_{j-i \ge {n-k+1}} (x_j^2 - x_i^2)}  
\det\left[\Delta M_{2n-2k+2}(\pm x|_j^{j+n-k};y_i)\right]_{i,j=1}^k  v_k(dY),
	$$
	where $\Delta$ is the following differential operator 
	$$
	    \Delta = - y\frac{d}{dy} + \varkappa(n,k),
	$$
	and $c(n,k), \varkappa(n,k)$ are the constants that depend on $n$ and $k$.
	
	If $G_n = \mathrm{SO}(2n+1)$ or  $G_n = \mathrm{Sp}(2n)$, then $\varkappa(n,k) = 0$ and
	$$
	c(n,k) = \frac{(2n-2k)!!}{(2n)!!} \prod_{i=0}^{k-1} {2n-2k+2i+2 \choose 2i+1}.
	$$
	
	In the case when $G_n = \mathrm{O}(2n)$ we have  $\varkappa(n,k) = 2(n-k)$ and 
	$$
	c(n,k) = \frac{(2n-2k-1)!!}{(2n-1)!!} \prod_{i=0}^{k-1} {2n-2k+2i+1 \choose 2i}.
	$$
}

\medskip

{\bf Remark.} Note that the derivative of the B-spline $M_m, m\in \mathbb{N}$ can be expressed as the difference of two splines of order $m-1$; namely, for any $m$ points $t_1,...,t_m \in \mathbb{R}$ the following equality holds (see \cite{Phillips}):
$$
	\frac{t_m-t_1}{m-1}\cdot \frac{d}{dt}M_m(t_1,...,t_m;t) = M_{m-1}(t_1,...,t_{m-1};t) - M_{m-1}(t_2,...,t_m;t).
$$

\medskip

{\bf 3. Proof of Theorem 1.}

{\it 3.1. The Laplace transform of orbital measures.}
    We define the Laplace transform of the orbital measure $\mu$ on the Lie algebra $\mathfrak{g}_n$ of the group $G_n$ as the orbital integral
    $$
	\widehat{\mu_X}(T) = \int_{\mathrm{Orbit}(X)} e^{\mathrm{tr} (TX)} \mu(dX) = 
	 \int_{G_n}  e^{\mathrm{tr}(T \, \mathrm{Ad}_g(X))} \, dg, 
$$
where $dg$ is the Haar measure on $G_n$.

    In general, the Laplace transform for the measures with a compact support is defined on the complexification of the Lie algebra. But let us notice that the funciton $\widehat{\mu_X}(T)$ is invariant under the adjoint action of $G_n$. Thus we can think of $T$ as of the matrix of canonical form and consider the Laplace transform of the orbital measure as a function on the coordinate space $\mathbb C^n$.
     
    Due to the Harish-Chandra theorem (see \cite{Harish-Chandra}, Theorem 2), one can obtain the following formulas for the Laplace transform of the orbital measure $\mu_X$ for $G_n$: 
    \begin{align}\label{fourier-abcd}
\widehat{\mu_X}^{(B_n, C_n)}(t_1,...,t_n) & =  \frac{(2n-1)!(2n-3)!...1!}{V_n(X^2)V_n(T^2)}\det\left[\frac{\sinh(t_ix_j)}{t_ix_j}\right]_{i,j=1}^n,  \\
\widehat{\mu_X}^{(D_n)}(t_1,...,t_n)& = \frac{(2n-2)!...2!}{V_n(X^2)V_n(T^2)} \det\left[\cosh(t_ix_j)\right]_{i,j=1}^n. 
\end{align}

    In the formulas \eqref{fourier-abcd}, (2) $B_n, C_n$ and $D_n$ are the standard notations for the classical series of Lie algebras of groups
    $\mathrm{SO}(2n+1)$, $\mathrm{Sp}(2n)$ and $\mathrm{O}(2n)$ respectively.
    
    These formulas imply that, up to multiplicative constants, the Laplace transform of the orbital measure $\mu_X$ has the form
    $$
	D_n(f; T; X) = \frac{\det\left[f(t_ix_j)\right]_{i,j=1}^n}{V_n(T^2)V_n(X^2)},	
    $$
    where $f$ is an analytic function with the Taylor series
    $$
        f(z) = \sum_{m=0}^{\infty} \frac{z^{2m}}{2^{2m}m!(\alpha+1)_m}.
    $$
    Here $(\alpha+1)_m = (\alpha+1)(\alpha+2)...(\alpha +m)$ is the Pochhammer symbol, and the parameter $\alpha$ is equal to $1/2$ in the case of the groups $\mathrm{SO}(2n+1)$ and $\mathrm{Sp}(2n)$ and $\alpha = -1/2$ when $G_n
= \mathrm{O}(2n)$.

    {\it 3.2. Projections and the Laplace transform.}
    For any Borel measure $\mu$ on $\mathfrak g_n$ the restriction of the Laplace transform to $\mathfrak g_k$ is equal to the Laplace transform of the measure $\mu^{(k0}$ which is the image of $\mu$ under projection onto $\mathfrak g_k$. Thus the problem is reduced to the computation of the quantities $D_n(f;t_1,...,t_k,0,...,0; X)$.
    
    {\bf Lemma 1.}
    {\it If $f$ is an even analytic function in some neighbourhood of 0, then the quantity $D_n(f;t_1,...,t_k,0,...,0; X)$ can be expressed in terms of divided differences of the functions $\varphi_i(y) = f(t_i\sqrt{y})$ so that 
    $$
	D_n(f; t_1,...,t_k, 0,...,0; X) = \frac{a_{n-k}(f)}{V_n(X^2)V_k(T^2)(t_1...t_k)^{2(n-k)}}  \times 
    $$ 

    \begin{equation}\label{divided}
	\times \prod\limits_{1\le j-i\le n-k} (x_j^2-x_i^2) \cdot \det\left[\varphi_i[x_j^2,...,x^2_{j+n-k}]\right]_{i,j=1}^k.
	\end{equation}
    
	Here $a_n(f) = \prod_{j=0}^{n-1} c_{2j}$, where $c_{2j}$ are the even coefficients of the Taylor series of $f$. 
    }
    
    The proof of this lemma is similar to the ones of Theorem 4.1 and 5.3 in the paper \cite{Faraut}.
    $\square$
    
    {\it 3.3. Doubling the knots.}
    Since the functions $\varphi_i$ from Lemma 1 depend on $\sqrt{z}$, the formula \eqref{divided} becomes inconvenient for applying the Hermite-Gennochi formula. But we can overcome this difficulty, if we double the the number of knots.
    
    {\bf Lemma 2.} {\it 
    Let $f(z)$ be an analytic function in a certain neighbourhood of $0$, and let $\varphi(z) = f(\sqrt{z})$. Then, given $m$ points $0< z_1 < ... < z_m$, the following equality for the divided differences holds:
    $$
		\varphi[z_1^2,...,z_m^2] = g[-z_m,...,-z_1,z_1,...,z_m],
	$$
	where  $g(z) = zf(z)$.
    }
    
    {\bf Proof.} By the definition of divided differences, we have
    $$
	\varphi[z_1^2,...,z_m^2] = \sum\limits_{l=1}^{m} f(z_l) \prod\limits_{r\ne l} \frac{1}{z_l^2 - z_r^2} = 
$$
$$ = 2 \sum\limits_{l=1}^{m} \frac{z_lf(z_l)}{2z_l} \prod\limits_{r\ne l} \frac{1}{z_l^2 - z_r^2} =g[-z_m,...,-z_1,z_1,...,z_m],
$$
 where $g(z) = zf(z)$. Lemma 2 is proved. $\square$
 
{\it 3.4. Applying the Hermite-Genocchi formula.}

    Due to Lemmas 1 and 2 the computation of $D_n(f; t_1,...,t_k, 0,...,0; X)$
has been reduced to the computation of the determinant of divided differences of the functions $g_i(z) = z f(t_iz)$. 

Applying the Hermite-Genocchi formula, we get the following equality:
\begin{multline}
	 g_i[-x_{j+n-k},...,-x_j,x_j,...,x_{j+n-k}] = \frac{1}{(2n-2k+1)!}\\
	\, \, \, \times \int\limits_{\mathbb{R}} (g_i(z))^{(2n-2k+1)} M_{2n-2k+2}(-x_{j+n-k},...,-x_j,x_j,...,x_{j+n-k};z) dz =\\
	= \frac{-1}{(2n-2k+1)!} \int\limits_{\mathbb{R}} (g_i(z))^{(2n-2k)} M'_{2n-2k+2}(\pm x|_j^{j+n-k};z) dz. \nonumber
\end{multline}

Now we want to express  $(g_i(z))^{(2n-2k)} = (zf_i(z))^{(2n-2k)}$ in terms of $f_i(z)=f(t_i z)$.

By the definition of $f$, for any $m \in \mathbb{N}\cup \{0\}$ we have
$$(g_i(z))^{(2n-2k)} =  
\sum_{m=0}^\infty \frac{t_i^{2m}z^{2m+1}}{2^{2m}m!(\alpha+1)_m} \cdot \gamma_{\alpha,m},$$
where
$$
	\, \gamma_{\alpha,m} = \gamma_\alpha (m) := \frac{(m+\frac{3}{2})_{n-k}}{(m+\alpha+1)_{n-k}}. 
	$$

The problem is now divided into two cases: \par 	
{\bf (B-C)} $G_n = \mathrm{Sp}(2n)$ or $\mathrm{SO}(2n+1)$.	In this case $\alpha = 1/2$, $\gamma_{\alpha,m} = 1$ and
$$
    (g_i(z))^{(2n-2k)} = t_i^{2n-2k} zf(t_iz).
$$
\par
{\bf (D)} $G_n = \mathrm{SO}(2n)$. In this case 
	$\alpha = -1/2$, 
	$\gamma_{\alpha,m} = 1 + \frac{2(n-k)}{2m+1}$
	and therefore,
$$
	(g_i(z))^{(2n-2k)} = t_i^{2n-2k} \cdot \left(zf(t_iz) + (2n-2k)t_i^{-1} f^{(-1)}(t_iz)\right),
$$
where $f^{(-1)}(\cdot)$ is the primitive of $f$ such that $f^{(-1)}(0)=0$.

	In both cases the divided differences of $g_i$'s can be expressed as follows:
	\begin{multline}\label{AfterHG}
 	g_i[-x_{j+n-k},...,-x_j,x_j,...,x_{j+n-k}] = \\ =\frac{t_i^{2n-2k}}{(2n-2k+1)!}
	\cdot \int\limits_{\mathbb{R}} f(t_iy) \, \Delta M_{2n-2k+2}(\pm x|_j^{j+n-k};y) dy,
\end{multline}
    where the operator $\Delta$ has the form
    $$
	\Delta = \Delta_{1/2} = -y\frac{d}{dy} \, \text{\, (the (B-C) case)}
	$$
	or
	$$
	\Delta = \Delta_{-1/2} = -y\frac{d}{dy} + 2(n-k) \, \, \text{\, (the (D) case)}.
	$$
    
    Now let us apply the Binet-Cauchy identity to the formula \eqref{AfterHG}:
    \begin{equation}\label{Binet1}
\begin{aligned}
    \det\left[  g_i[-x_j,...,-x_{j+n-k},x_{j+n-k},...,x_{j}]\right] =  \prod_{i=1}^k t_i^{2n-2k-1} \cdot \frac{1}{((2n-2k+1)!)^k} \times
\\
\\ \times \int\limits_{\mathbb{R}^k_+}
\det\left[f(t_i y_j))\right]_{i,j=1}^k \det\left[\Delta M_{2n-2k+2}(\pm x|_j^{j+n-k}; y_i)\right]_{i,j=1}^k dy_1...dy_k.
\end{aligned}
\end{equation}

    To complete the proof of Theorem 1 it suffices to compare the LHS and RHS in \eqref{divided} and to apply the inverse Laplace transform to both sides of that equation.  Theorem 1 is proved. $\square$.
    
\smallskip

{\bf Remark.} After this work was finished, the author became aware of an interesting paper of Defosseux \cite{D}, where she considered, for the matrix $X \in \mathfrak g_n$, the array $M(X)$ of the union of spectra of $X$ and of all its images under projections $p_k^n, k=1,...,n-1$. In turns out that, the radial part of projection of the orbital measure $\mu_X$ can be derived as a correlation function of the determinantal point process on arrays $M(X)$ with the kernel given explicitely (see Theorem 6.3 of the work of Defosseux).

\medskip

{\bf Acknowledgement.} I am deeply grateful to Grigori Olshanski, who suggested me this problem, for many useful discussions and paying attention to this work.\\[2mm]
    
{Steklov Mathematical Institute of the Russian Academy of Sciences,\\
8, Gubkina str., Moscow, 119991, Russian Federation}\\

{\bf e-mail}: dmitry.zubov.93@gmail.com

\end{document}